\documentclass[a4paper,fleqn]{cas-dc}

\usepackage[authoryear,longnamesfirst]{natbib}
\usepackage{amsmath}
\usepackage{amssymb}
\usepackage{bm}
\usepackage{graphicx}
\usepackage{subfigure}

\usepackage{bm}
\usepackage{amsmath,amssymb,amsfonts}

\usepackage[utf8]{inputenc}
\usepackage[T1]{fontenc}
\usepackage{etoolbox}

\usepackage{algorithmic}
\usepackage{graphicx}
\usepackage{textcomp}

\def\tsc#1{\csdef{#1}{\textsc{\lowercase{#1}}\xspace}}
\tsc{WGM}
\tsc{QE}
\tsc{EP}
\tsc{PMS}
\tsc{BEC}
\tsc{DE}

\begin{document}
\let\WriteBookmarks\relax
\def\floatpagepagefraction{1}
\def\textpagefraction{.001}

\shorttitle{Finite-time stochastic control}

\shortauthors{Xiaoxiao Peng et~al.}

\title [mode = title]{Finite-time stochastic control for complex dynamical systems: The estimate for control time and energy consumption}                      


%
\author[1,2]{Xiaoxiao Peng}[
                        style=chinese         
                        ]



\ead{xxpeng19@fudan.edu.cn}



\affiliation[1]{organization={Shanghai Center for Mathematical Sciences, Fudan University},
    addressline={2005 Songhu Road}, 
    city={Shanghai},
    postcode={200438}, 
    country={China}}

\affiliation[2]{organization={School of Mathematical Sciences, Fudan University},
    addressline={220 Handan Road}, 
    city={Shanghai},
    postcode={200433}, 
    country={China}}

\affiliation[3]{organization={Research Institute of Intelligent Complex Systems, Fudan University},
    addressline={220 Handan Road}, 
    city={Shanghai},
    postcode={200433}, 
    country={China}}


\author[3]{Shijie Zhou}[style=chinese]

\cormark[1]

\fnmark[1]
\ead{szhou14@fudan.edu.cn}




\cortext[cor1]{Corresponding author}



\begin{abstract}
Controlling complex dynamical systems has garnered significant attention within academic circles in recent decades. While most existing works have focused on closed-loop control schemes with infinite-time durations, this paper introduces a novel finite-time, closed-loop stochastic controller that considers control time, energy consumption, and their dependence on system parameters. This stochastic control technique not only enables finite-time control in chaotic dynamical systems but also facilitates finite-time synchronization in unidirectionally coupled systems. Importantly, our proposed scheme offers several advantages over existing deterministic finite-time controllers, particularly in terms of time and energy consumption. Through numerical experiments utilizing random ecosystems, neural networks, and Lorenz systems, we provide evidence of the effectiveness of our analytical results. We anticipate that this stochastic scheme will find wide-ranging applications in the control of complex dynamical systems and in achieving network synchronization.\end{abstract}



\begin{keywords}
finite-time, stochastic control, control time and energy consumption
\end{keywords}

\maketitle
\newtheorem{theorem}{Theorem}
\newtheorem{remark}{Remark}
\newtheorem{example}{Example}
\newtheorem{condition}{Condition}
\newtheorem{corllary}{Corllary}
\newtheorem{lemma}{Lemma}
\newtheorem{definition}{Definition}
\section{Introduction}
In recent years,  controlling complex dynamical systems has attracted much interest. Among various techniques developed for controlling, the methodologies of open-loop have received  much attention (refer to \cite{b11, b25, b26, b27, b28, b33, b34}). In the open-loop control, the agent selects several nodes in advance and applies predefined state-independent control signals to drive the system from an initial state to a desired target (refer to \cite{b45,b46,b47}). It is direct and can reduce the complexity in many cases.  However, such controls are usually energy-consuming, unsustainable, and less robust. To overcome these difficulties, the closed-loop control methods are proposed and applied in many complex dynamical systems. To be precise, several feedback schemes of close-loop control, have been validated to be effective numerically or theoretically in controlling various kinds of systems with or without time delay/noise (refer to \cite{b2,b3,b8,b4,b5,b6,b7,b9,b10,b58,b61,b62}).

Most closed-loop control schemes are often assumed to be deterministic, overlooking the presence of noise in real-world systems. Therefore, designing closed-loop control schemes for systems that are affected by noise is imperative. Numerous analytical and numerical investigations have shown that closed-loop control schemes can facilitate stochastic stabilization or synchronization of noisy systems (see \cite{b4,b5,b6,b7,b12,b13,b14,b39, b40,b60,b71,b72}). However, current closed-loop control schemes have limited control flexibility and require infinite control time. A few non-closed stochastic control schemes have been developed that can control specific systems within a finite time frame (see \cite{b42, b43, b44}). Hence, it is natural to ask whether we can design a finite-time, closed-loop stochastic controller for general systems.

Numerous researchers have made significant strides in exploring finite-time stabilization of stochastic dynamical systems. For instance, some studies have established finite-time criteria for switched nonlinear stochastic systems, as evidenced in \cite{b49,b55,b56}. Others have devised Lyapunov criteria for finite-time stability of stochastic systems, as exemplified in \cite{b50,b52,b54}. Furthermore, certain studies have investigated stochastic controllers utilizing backstepping design to realize finite-time stability, as elaborated in \cite{b51, b54}. However, to our best knowledge, there remains a lack of research examining the estimate for control time and energy consumption for finite-time stochastic controllers applied to general dynamical systems.



In this paper, the above questions are fully addressed. Drawing inspiration from finite-time closed-loop deterministic control schemes, as documented in prior literature such as \cite{b15,b16,b17,b18,b19,b20,b21,b59,b57}, we propose a stochastic scheme that can be employed to stabilize complex dynamical systems with high probability in a finite time frame. Notably, our work differs from past literature (as highlighted in \cite{b20,b15,b59}) in that the controller is added on the ${\rm d}B_t$ term rather than the ${\rm d}t$ term, which makes the analytical estimation more difficult (See Remark \ref{remark2}). The impact of the coupling gain and parameters on control time and energy is studied using several real biophysical systems, including random ecosystems, neural networks and Lorenz systems. We ascertain upper bounds for control time and energy consumption via analytical means. Our scheme aims to achieve not only finite-time control in chaotic dynamical systems, but also finite-time synchronization in unidirectionally coupled systems. Intriguingly, compared to deterministic finite-time controllers, our new scheme depicts notable advantages from a physical perspective concerning time and energy consumption.


The rest of this paper is structured as follows. In Section \ref{pre}, we introduce the stochastic model that we seek to modulate and the closed-loop control schemes. In Section \ref{example}, we establish upper bounds for both control time and energy consumption, and we utilize several illustrative examples to demonstrate the applicability of our schemes in stabilization and synchronization. In Section \ref{comparison}, we compare the proposed stochastic controller with the corresponding deterministic controller in terms of time and energy consumption. In Section \ref{proof}, we provide mathematically rigorous proofs to support the conclusions obtained in Section \ref{example}. Finally, in Section \ref{conclusion}, we offer concluding remarks and potential directions for future research.

\textbf{Notations}. For any $\bm{x},\bm{y}\in\mathbb{R}^n$, $\|\bm{x}\|$ represents the standard Euclidean norm of $\bm{x}$, while $\langle\bm{x},\bm{y}\rangle$ refers to the standard Euclidean inner product of $\bm{x}$ and $\bm{y}$. ${B}_{t}( t\geq 0)$ is a one-dimensional Brownian motion. Moreover, we denote by $(\Omega,\mathscr{F},\{\mathscr{F}_t\}_{t\geq 0},\mathbb{P})$ the complete probability space with a filtration $\mathscr{F}_t$ satisfying the usual condition (${\rm i.e.}$, it is increasing and right continuous while $\mathscr{F}_0$ contains all $\mathbb{P}-$null sets) where the one-dimensional Brownian motion $B_t$ is defined (See Eq. \eqref{1}). For a random variable $\xi$, $\mathbb{E}\xi$ denotes its expectation. For a measurable set $\bm{A}$, $\bm{1}_{\bm{A}}$ denotes the indicator function of $\bm{A}$. ${\rm a.s. }$ represents almost surely in abbreviation. For $x\in\mathbb{R}$, ${\rm Sgn}(x)$ denotes the sign of $x$, which
 indicates ${\rm Sgn}(x) = 1,0$ and $-1$ for $x > 0, x = 0$ and $x < 0$, respectively.

\section{Model Formulation}\label{pre}

We consider a $n$-dimensional stochastic controlled dynamical systems as follows.
\begin{equation}\label{1}
	{\rm d}\bm{x}_t=\bm{f}(\bm{x}_t){\rm d}t+\bm{u}(\bm{x}_t){\rm d}B_t
\end{equation}
with initial state $\bm{x}(0)=\bm{x}_0$ is a random variable, where $\bm{x}_t\in\mathbb{R}^n$ is the state vector, $\bm{f}(\bm{x})$ denotes the vector field and $\bm{u}(\bm{x})$ denotes the stochastic controller to be designed. ${B}_t$ is a one-dimensional Brownian motion and the filtration $\mathscr{F}_t(t\geq 0)$ is generated by $\bm{x}_0$ and the Brownian motion $B_s(s\leq t)$.
According to the classical theory of stochastic differential equations, we need locally Lipschitizian conditions for vector $\bm{f}$ as follows.
\begin{condition}\label{c1}
	(Locally Lipschitzian condition)
	For any given positive number $m$, there exists a positive number $K_m$ such that 
	$$
	\|\bm{f}(\bm{x})-\bm{f}(\bm{y})\|\leq K_m\|\bm{x}-\bm{y}\|
	$$
	for any $\|\bm{x}\|\leq m$ and $\|\bm{y}\|\leq m$. $\|\cdot\|$ represents for Euclidean norm in $\mathbb{R}^n$.
\end{condition}
We suppose that $\bm{x}_t\equiv \bm{0}$ is a solution for origin system. We also need conditions as follows, which has been proposed and investigated in past literatures (see \cite{b35, b36, b37, b38}).
\begin{condition}\label{c2}
	(Globally one-sided Lipschitzian condition)
	There exists a positive number $L$ such that
	$$
	\langle \bm{x},\bm{f}(\bm{x})\rangle \leq L \|\bm{x}\|^2
	$$
	for any $\bm{x}\in\mathbb{R}^n$.
\end{condition}

Denote by the  stopping time  when the trajectory first hits zero.
$$
\tau\triangleq \inf \bigg\{t>0\bigg| \bm{x}_t = \bm{0}\bigg\}.
$$
\begin{definition}\label{def1}
	Eq. \eqref{1} is said to be finite-time stable if $\mathbb{P}(\tau<+\infty)=1$. Equivalently speaking, $\tau$ exists finitely in the sense of probability one.
\end{definition}
\begin{definition}
	If Eq. \eqref{1} is finite-time stable,  for $q>0$, denote by the $L^q$ energy consumption
	$$
	\mathcal{E}_q\triangleq\int_0^\tau \|\bm{u}(\bm{x}_s)\|^q{\rm d}s.
	$$
 If $\mathbb{E}\mathcal{E}_q<+\infty$, we say that the energy of the system \eqref{1} is $L_q$ integrable.
\end{definition}
\begin{remark}
In  framework of deterministic control, the $L_q$ energy is usually defined for $q>1$. However, in our framework of stochastic control, we usually consider the situation $q<1$ (See Theorem \ref{theorem2}).
\end{remark}
In most existing closed-loop stochastic control systems, the feedback controller has been designed as $\bm{u}(\bm{x})=k\bm{x}$, where the coupling gain $k$ satisfies the condition $k > \sqrt{2L}$ (refer to \cite{b48,b4}). This type of stochastic controller can asymptotically exponentially stabilize the solution of \eqref{1} to the zero state under Condition \ref{c2}. However, the convergence time is infinite. Therefore, we aim to find a method to achieve a finite control time. As has been done in deterministic dynamical systems (refer to \cite{b15}), we propose the following controller:
\begin{equation}\label{2}
	\bm{u}(\bm{x})=k\bm{x}\cdot\bm{1}_{\|\bm{x}\|\geq 1}+ k\|\bm{x}\|^{\alpha-1}\bm{x}\cdot\bm{1}_{\|\bm{x}\|< 1},
\end{equation}
where $\alpha\in (0,1)$ denotes the steepness exponent, and $k>\sqrt{2L}$ is the coupling gain.

\begin{figure}
	\includegraphics[width=0.45\textwidth]{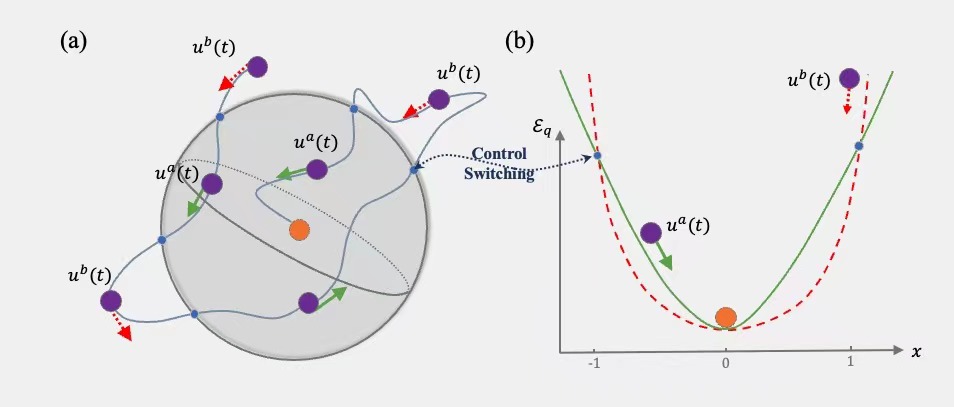}
	\caption{Physical underpinning of $\bm u(\bm x)$ for the three- and one-dimension case, respectively. Here, $\bm u^{a}(\bm{x})=k\bm x \Vert\bm x\Vert^{\alpha-1}$, $\bm u^{b}(\bm{x})=k\bm x$.
	Since controller \eqref{2} is designed as a stochastic controller, it is possible for the trajectory to enter and exit the unit ball multiple times before eventually converging to zero within a finite time (refer to Remark \ref{remark2}). This behavior is significantly different from the deterministic controller discussed in \cite{b15}.}\label{physical}
\end{figure}

By replacing the diffusion term $\bm u(\bm x_{t}){\rm d}\bm B_{t}$ with $-\bm u(\bm x_{t}){\rm d}t$, the equation~\eqref{1} can be transformed into a closed-loop control for a deterministic system, which has been proven to effectively stabilize the original ordinary differential equations ${\rm d}\bm{x}_t=\bm{f}(\bm{x}_t){\rm d}t$ within a finite time duration. \textcolor{black}{It is important to note that this deterministic controller is similar to, but distinct from, the one established in \cite{b15} (refer to Remark \ref{remarkjia2}).} The physical basis of $\bm u(\bm x)$ is illustrated in Figure~\ref{physical}.  Moreover, it can be applied to achieve synchronization. But to our best knowledge, there have not been systematic investigations on the time and energy consumption for the stochastic closed-loop controller~\eqref{2}.

\textcolor{black}{\begin{remark}\label{remarkjia3}
Since Brownian motion $\bm{B}_t$ is symmetric about zero, the control $\bm{u}(\bm{x}_t){\rm d}\bm{B}_t$ is equivalent to $-\bm{u}(\bm{x}_t){\rm d}\bm{B}_t$. This represents a significant distinction between the stochastic controller and the deterministic controller $-\bm{u}(\bm{x}_t){\rm d}t$.
\end{remark}}

\textcolor{black}{\begin{remark}\label{remarkjia2}
In \cite{b15}, the authors established the following finite-time closed-loop control
$$
\dot{\bm{x}}=\bm{f}(\bm{x})-\bm{v}(\bm{x}),
$$
where the controller $\bm{v}(\bm{x})=[{v}_1(\bm{x}),{v}_2(\bm{x}),\cdots,{v}_n(\bm{x})]^{\rm T}$ is denoted by
$$
{v}_i(\bm{x})\triangleq\left\{ \begin{aligned}&kx_i, \|\bm{x}\|>1,\\&k ({\rm sgn}~x_i) |x_i|^\alpha, \|\bm{x}\|\leq 1.    \end{aligned}\right.
$$
It is worthwhile to note that $\bm{u}(\bm{x})$ is similar to, but distinct from $\bm{v}(\bm{x})$. Additionally, it should be emphasized that $\bm{u}(\bm{x})$ is continuous in $\mathbb{R}^n$, whereas $\bm{v}(\bm{x})$ is not continuous on the surface of the unit ball $\|\bm{x}\|=1$. The reason why we choose $\bm{u}(\bm{x})$ instead of $\bm{v}(\bm{x})$ lies in that the major mathematical tool used in the proof of Theorem \ref{theorem1} and \ref{theorem2} is the generalized It$\hat{o}$'s formula (Lemma \ref{lemma1}), which requires the vector fields to be continuous.
\end{remark}}

Here, we first show that it can achieve stabilization as well as synchronization via some examples and then give the detailed proofs.

\section{Realization of Stabilization and Synchronization}\label{example}
\subsection{Stabilization}

To begin with, we will introduce two theorems that facilitate our characterization of convergence time and energy consumption estimates. The proof of the theorems will be provided in Section \ref{proof}.

\begin{theorem}\label{theorem1}
	Assume that Conditions \ref{c1},\ref{c2} are valid. Then Eq. \eqref{1} with stochastic controller \eqref{2} is finite-time stable for $k>\sqrt{2L}$. Furthermore, the time estimation upper bound $T_f^{\rm Sup}$ satisfies that

$$
\mathbb{E}\tau\leq T_{f}^{\rm Sup}\triangleq\left\{
\begin{aligned}
	&\dfrac{2\mathbb{E}(\log\|\bm{x}_0\|;\bm{1}_{\|\bm{x}_0\|>1})}{{k^2}-2L}\\
 &+\dfrac{8k^2}{(k^2-2L)^2}[\mathbb{P}(\|\bm{x}_0\|>1)\\
&+\mathbb{E}(\|\bm{x}_0\|^{p^*};\bm{1}_{\|\bm{x}_0\|\leq 1})],
 \alpha<\dfrac{3}{4}+\dfrac{L}{2k^2},\\
	&\dfrac{2\mathbb{E}(\log\|\bm{x}_0\|;\bm{1}_{\|\bm{x}_0\|>1})}{{k^2}-2L}\\
 &+\dfrac{[\mathbb{P}(\|\bm{x}_0\|>1)
+\mathbb{E}(\|\bm{x}_0\|^{2-2\alpha};\bm{1}_{\|\bm{x}_0\|\leq 1})]}{(1-\alpha)[(2\alpha-1)k^2-2L]},\\
	&\alpha\geq \dfrac{3}{4}+\dfrac{L}{2k^2},
\end{aligned}\right.
$$
where $p^*\triangleq\dfrac{1}{2}-\dfrac{L}{k^2}$. 
\end{theorem}

\begin{theorem}\label{theorem2}
	Assume that Conditions \ref{c1},\ref{c2} are valid. The system \eqref{1} with the controller \eqref{2} is finite-time stable. When  $0<q<\min\left\{2-2\alpha,1-\dfrac{2L}{k^2}\right\}$,  the energy is $L_q$-integrable. Furthermore, we have the estimation of the $L_q$ energy consumption as follows.
	$$
 \begin{aligned}
\mathbb{E}\mathcal{E}_q\leq\mathcal{E}_q^{\rm Sup}&\triangleq-\dfrac{k^q}{H_2 (q)}[\mathbb{E}(\|\bm x_0\|^q )+2\mathbb{P}(\|\bm x_0\|>1)\\
&+\mathbb{E}(\|\bm x_0\|^q ;\|\bm x_0\|\leq 1)],
 \end{aligned}
	$$
where $H_2(q)\triangleq qL+\dfrac{q(q-1)k^2}{2}<0.$
 
\end{theorem}

It can be easily checked  $T_{f}^{\rm Sup}$ and $\mathcal{E}_q^{\rm Sup}$ are decreasing with respect to the parameter $k$ when $\alpha$ and $\bm x_{0}$ are fixed. 
Let $k\to \sqrt{2L}+$, the time estimation $T_{f}^{\rm Sup}$ will diverge, which aligns with our intuition. For the scalar system ${\rm d}x_t=Lx_t+kx_t{\rm d}B_t$, it follows from It$\hat{o}$'s formula that ${\rm d}\log x_t=(L-\dfrac{k^2}{2}){\rm d}t+k{\rm d}B_t$, indicating that the system converges or diverges exponentially at a rate of $L-\dfrac{k^2}{2}$. As $k\to \sqrt{2L}+$, the system converges at a slower rate. Since  our controller $\eqref{2}$ aligns with $kx$ when $|x|>1$, it is natural for  the time estimation $T_{f}^{\rm Sup}$ to diverge. 

\begin{figure}
		\includegraphics[width=0.55\textwidth]{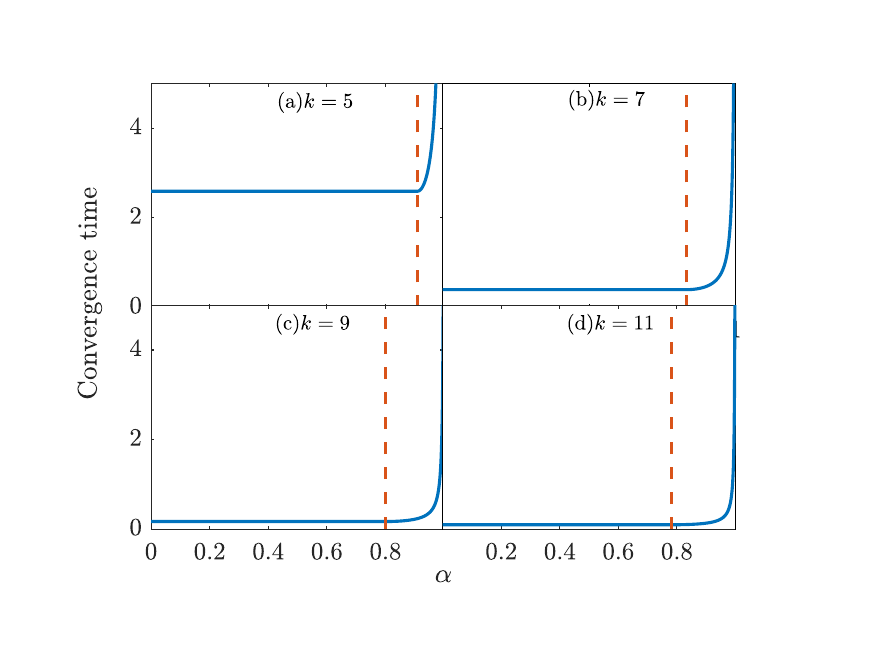}
		\caption{The solid (blue) cureves indicate the dependence of the $T_{f}^{\rm Sup}$ on $k$ and $\alpha$. And the dashed (red) curves indicate  $\alpha=3/4+L/(2k^{2})$. The parameters are set as $L=8$, $\Vert \bm x_{0} \Vert=3$.}\label{fig_J}
	\end{figure}

 We now discuss how $\alpha$ influences the  convergence time. We firstly note that $T_{f}^{\rm Sup}$ is continuous over $\alpha$.  When $\dfrac{3}{4}+\dfrac{L}{2k^2}\leq\alpha<1$, we first analyze the monotonicity of the function $J(\alpha)\triangleq(1-\alpha)[(2\alpha-1)k^2-2L]$. It can be calculated that 
$$
\frac{{\rm d}J}{{\rm d}\alpha}=-4k^2\alpha+3k^{2}+2L.
$$
So $\alpha=\dfrac{3}{4}+\dfrac{L}{2k^2}$ is the maximum point of function $J(\alpha)$. And  $\mathbb{E}(\|\bm{x}_0\|^{2-2\alpha};\bm{1}_{\|\bm{x}_0\|\leq 1})]$ is monotonically increasing over $\alpha \in (0,1)$ obviously,    which implies that $T_{f}^{\rm Sup}$ achieves minimum when $\alpha \in \left (0, \dfrac{3}{4}+\dfrac{L}{2k^2}\right]$. Addtionally, when $\alpha=1$, the controller specified in \eqref{2} becomes a infinite-time controller and this aligns with the fact that the time estimation $T_{f}^{\rm Sup}$ diverges when $\alpha\to 1-$ (See Figure. \ref{fig_J}). However, for $\alpha\in [0,\dfrac{3}{4}+\dfrac{L}{2k^2})$, $T_{f}^{\rm Sup}$ is a constant, which implies that the time estimation is conservative, leaving room for further improvement.

We give some illustrative examples in the following.

\begin{example}\label{examplenew1}
Consider a one-dimensional linear system described as
\begin{equation}\label{newexample}
{\rm d}x_t=Lx_t{\rm d}t+u(x_t){\rm d}B_t,
\end{equation}
where $L=2$. According to Theorem \ref{theorem1}, when $k>2$, the system \eqref{newexample} is finite-time stable. As shown in Figure. \ref{figure4a}, the convergence time and energy consumption show a decreasing tendency as $k$ increases when fix $\alpha$. As seen in Figure. \ref{figure4b}, the convergence time shows increasing tendency as $\alpha$ increases from $0$ to $1$ and diverges when $\alpha\to 1-$, while the energy consumption does not have an obvious monotonic trend when $\alpha$ varies between $(0,1)$. All these align with our estimations obtained in Theorem \ref{theorem1} and \ref{theorem2}.
\begin{figure}
\subfigure[]{\includegraphics[width=0.25\textwidth]{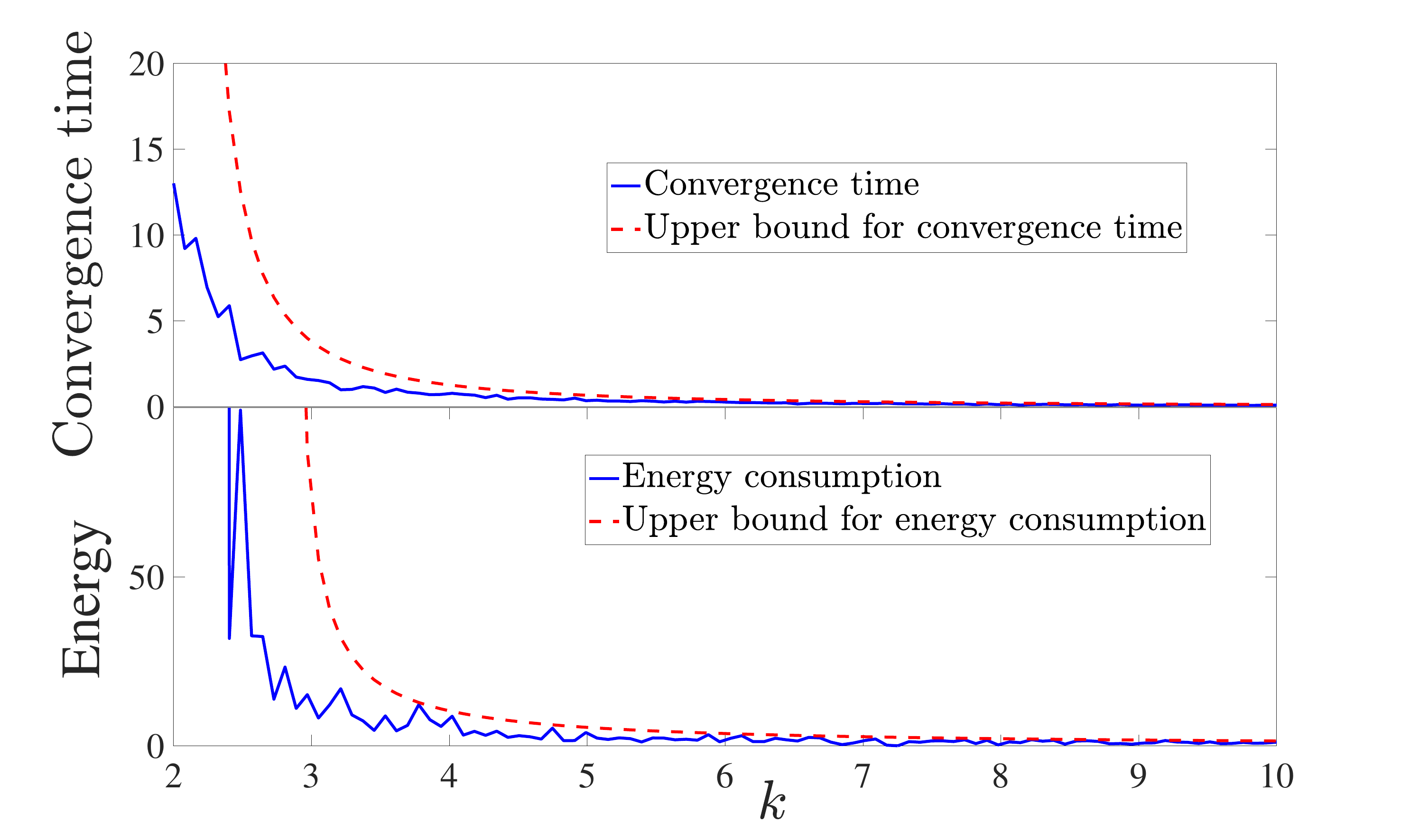}\label{figure4a}}\subfigure[]{\includegraphics[width=0.25\textwidth]{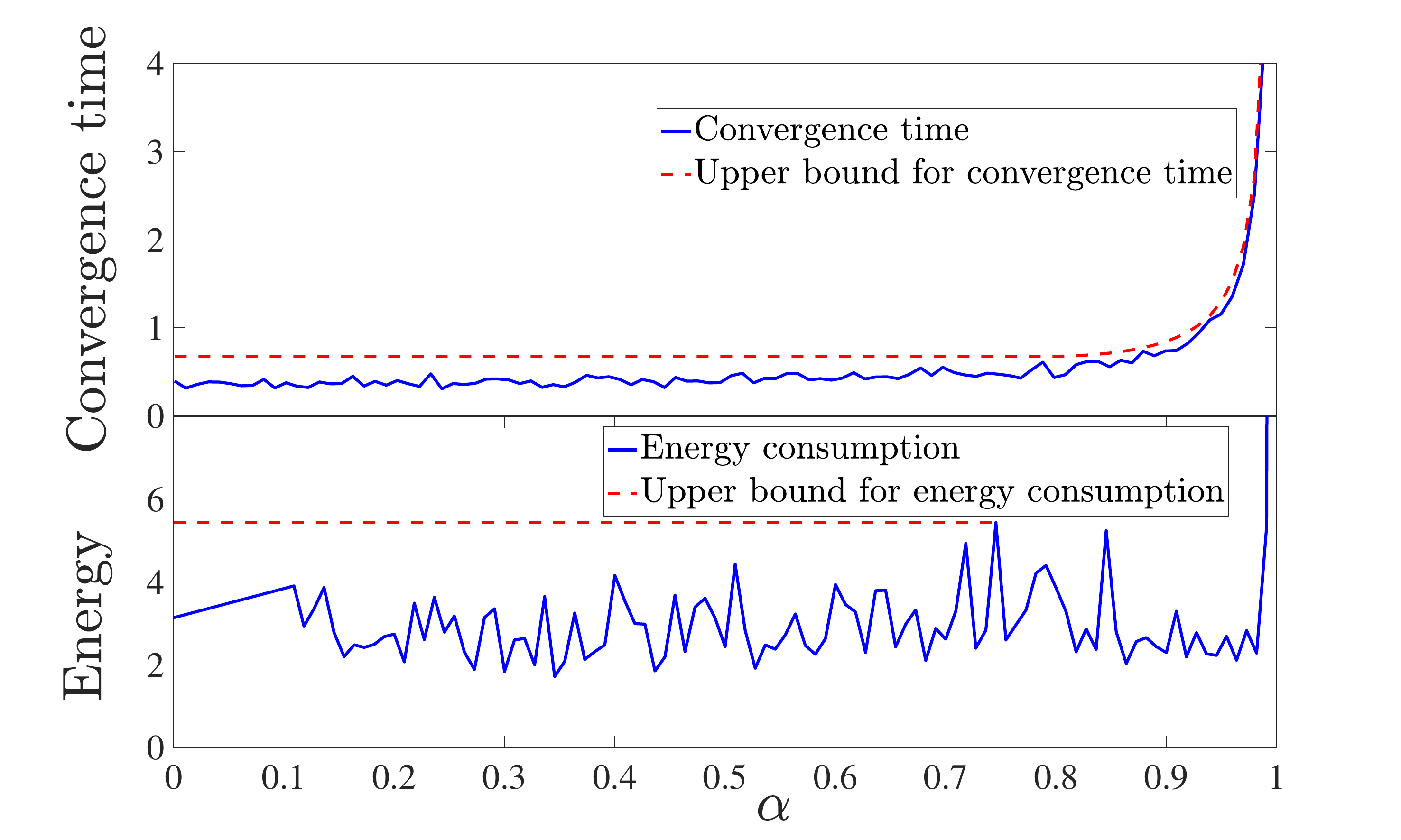}\label{figure4b}}
\caption{The convergence time and energy consumption for system \eqref{newexample} for different values of $\alpha$ and $k$. The solid (blue) curves represent the numerical evaluation of convergence time and energy consumption, which is counted as the average of \textcolor{black}{1000} random realization of system \eqref{newexample}. The dashed (red) curves indicate the upper bound obtained from Theorem \ref{theorem1} and \ref{theorem2}.
(a) When fix $\alpha$, convergence time and energy consumption shows a decreasing  tendency when $k$ increases. The parameters are $\alpha=0.5$, $q=1/2$. (b)When fix $k$, convergence time shows an increasing tendency when $\alpha$ increases while energy consumption does not have an obvious monotonic trend when $\alpha$ varies. The parameters are $k=5, q=1/2$.  For both (a) and (b), the initial value is $x_0=10$.  \textcolor{black}{Here, Euler-Maruyama scheme (see \cite{b29, b30}) with step size $\triangle t=10^{-6}$    is used for integrating the system of stochastic differentian equation \eqref{newexample}.}    }\label{figure4}
\end{figure}
\end{example}

\begin{example}\label{example22}
	Consider a two-dimensional nonlinear neural networks described as
 \begin{equation}\label{example1}
		{\rm d}\bm{x}_t=-\bm{C}\bm{x}_t{\rm d}t+\bm{A}\bm{g}(\bm{x}_t){\rm d}t+\bm{u}(\bm{x}_t){\rm d}B_t,
	\end{equation}
	where 
	$$
	\bm{C}=\left[\begin{array}{ll}
		1&2\\
		3&4
	\end{array}\right], \bm{A}=\left[\begin{array}{ll}
		3&3\\
		1&3
	\end{array}\right],$$ and the activation function $\bm{g}(\bm{x})=[{g}_1(\bm{x}_1),{g}_2(\bm{x}_2)]^{\rm T}$ satisfies that $g_1(x)=\tanh(x)$ and $g_2(x)=\tanh(2x)$.
	
	Using $|\tanh(x)|\leq |x|$ and $|\tanh(2x)|\leq 2|x|$, we obtain that
	$$\begin{aligned}
		&\langle \bm{x}, -\bm{C}\bm{x}+\bm{A}\bm{g}(\bm{x})\rangle
		\\ &=-x_1^2-4x_2^2-5x_1x_2+3x_1\tanh(x_1)\\&+3x_2\tanh(2x_2)+x_2\tanh(x_1)+3x_1\tanh(2x_2)
		\\  &\leq 2x_1^2+2x_2^2+12|x_1x_2|\leq 8(x_1^2+x_2^2).
	\end{aligned}
	$$
	Thus, the Lipschitzian constant value $L$, which is defined in Condition \ref{c2}, can be taken as $8$. According to Theorem \ref{theorem1}, we obtain that for each pair of $k>4$ and $\alpha\in(0,1),$ the controlled system \eqref{example1} is finite-time stable. As shown in Figure \ref{fig1}, the numerical evaluation of the expectation of convergence time, $\tau_k$, shows a decreasing tendency with the increase of $k$. The numerical average value, $\tau_k$, is always much smaller than the analytical upper bounds derived in Theorem \ref{theorem1}, indicating that estimation presented in our theorem is only in an approximating manner, leaving room for further improvement.
	
	\begin{figure}
		\includegraphics[width=0.45\textwidth]{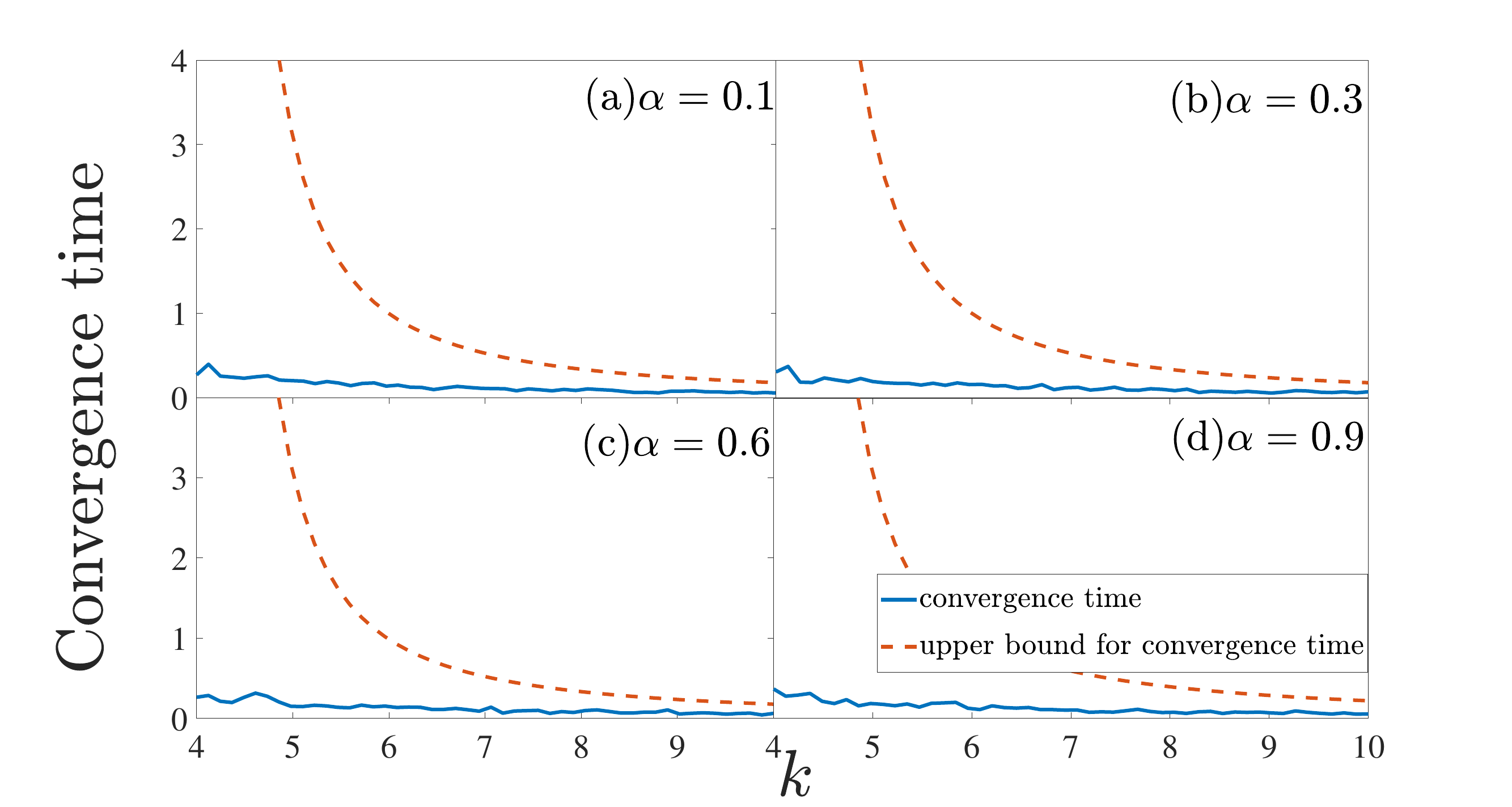}
		\caption{Dependence of the mean value of the convergence time $\tau_k$ on the coupling gain $k$ and $\alpha$. The solid (blue) curves indicate the numerical evaluation of $\tau_k$, which is counted as the average of 100 random realizations of system \eqref{example1}. The dashed (red) curves indicate upper bound obtained from Theorem \ref{theorem1}. The initial states are $\bm{x}_0=[10,10]^{\rm T}$.}\label{fig1}
	\end{figure}
	\end{example}

\begin{example}\label{example33}
	Consider May's classic ecosystem 
	\cite{b24, b41} described by 
	\begin{equation}\label{example2}
		{\rm d}\bm{x}_t=\bm{C}\bm{x}_t{\rm d}t+\bm{u}(\bm{x}_t){\rm d}B_t,
	\end{equation}
	where each species $x_i$ is one-dimensional, $\bm{C}=(c_{ij})_{N\times N}$ describes the random mutual interactions with $c_{ii}=-r$, and $N$ is the population size. The off-diagonal elements $c_{ij}$ are set as mutually independent Gaussian random variables $\mathscr{N}(0,\sigma^2)$ with probability $p$ and the probability for the elements to be zero is $1-p$.
	
	Denote each element of matrix $\bm{E}\triangleq\dfrac{1}{2}\left(\bm{C}^{\rm T}+\bm{C}\right)$ by $e_{ij}=\dfrac{1}{2}(c_{ij}+c_{ji})$. The expectation is $\mathbb{E}e_{ij}=0$ and the variance is given by 
	$$
	\mathbb{D}e_{ij}=\dfrac{1}{4}\mathbb{D}(e_{ij}+e_{ji})=\dfrac{1}{2}\mathbb{D}e_{ij}=\dfrac{1}{2}\mathbb{E}e^2_{ij}=\dfrac{1}{2}p\sigma^2.
	$$
	It follows from the semicircle law for random matrices that the eigenvalues of $\bm{E}$ are located in $[-r-\sqrt{2Np}\sigma,-r+\sqrt{2Np}\sigma]$ in a probabilistic sense as $N\to+\infty$.
	
	Since $\left\langle\bm{x},\bm{C}\bm{x}\right\rangle=\left\langle\bm{x},  \dfrac{1}{2}\left(\bm{C}^{\rm T}+\bm{C}\right)\bm{x}\right\rangle$, the Lipschitzian constant $L$ can be taken as $\eta_{\max}\triangleq-r+\sqrt{2Np}\sigma$ approximately. According to Theorem \ref{theorem1}, we obtain that when $k>\sqrt{2\eta_{\max}}$, the controlled system \eqref{example2} is finite-time stable.
	\begin{figure}
		\includegraphics[width=0.5\textwidth]{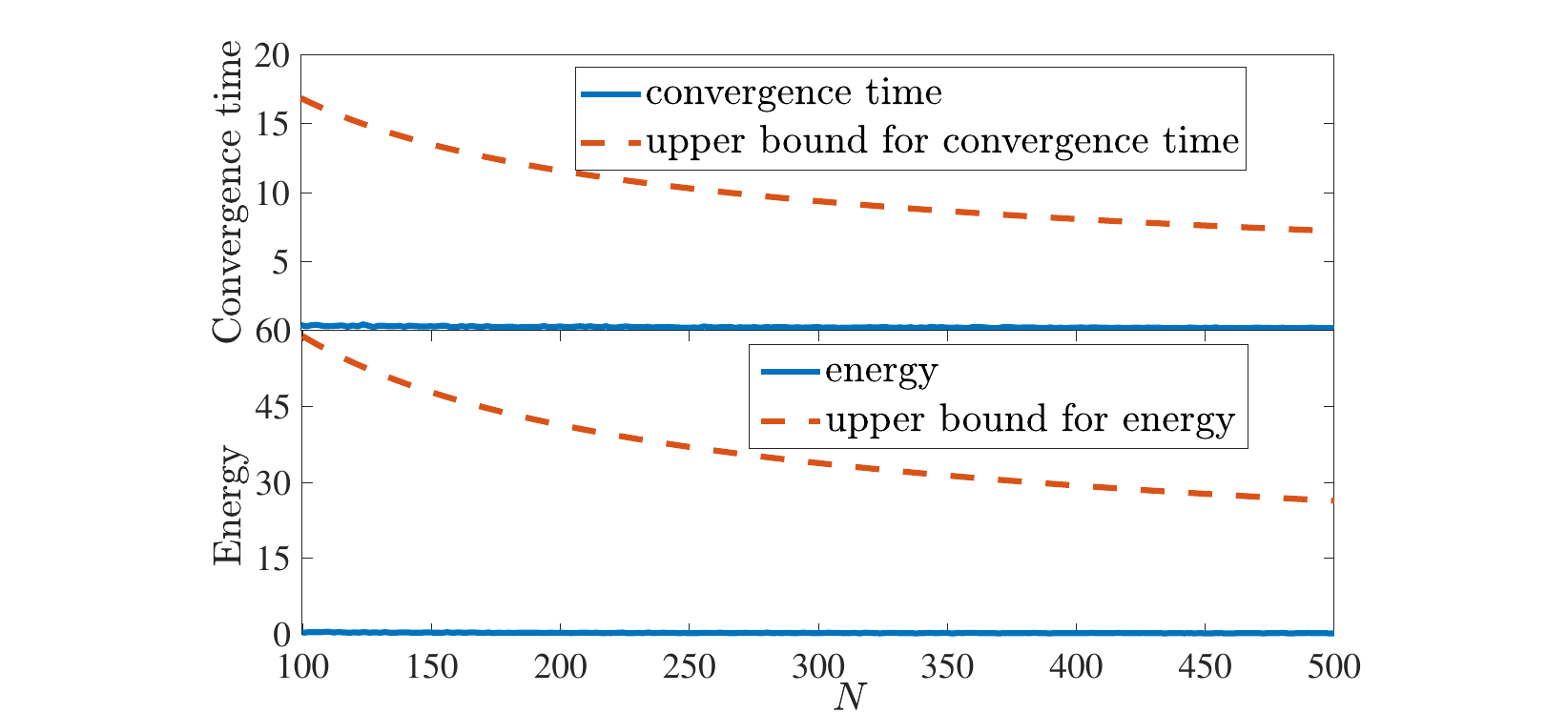}
		\caption{The convergence time and the required energy consumption for \eqref{example2}, vary respectively, with the increase of $N$. The parameters are $p=\dfrac{1}{3}$, $\sigma=1$, $r=1$, $q=0.1$, $k=1.1\sqrt{2\eta_{\max}}$, and the initial states are $\bm{x}_0=[1,1\cdots,1]^{\rm T}$. The solid (blue) curves indicate the numerical evaluation of $\tau$ and $\mathcal{E}_q$, which is counted as the average of 100 random realizations of system \eqref{example2}. The dashed (red) curves indicate upper bound obtained from Theorem \ref{theorem1} and \ref{theorem2}.}\label{fig2}
	\end{figure}
\end{example}
\begin{remark}
As shown in Figure. \ref{figure4}, our analytical upper bound slightly exceeds numerical results in Example \ref{examplenew1}.
 However, Figures \ref{fig1} and \ref{fig2} demonstrate that in Examples \ref{example22} and \ref{example33}, respectively, our analytical upper bound significantly surpasses the numerical results. The reason for this phenomenon lies in our use of the inequality $\langle\bm{x},\bm{f}(\bm{x})\rangle\leq L\|\bm{x}\|^2$ in obtaining the estimations in Theorems \ref{theorem1} and \ref{theorem2}.  However, this estimation is excessively conservative for Examples \ref{example22} and \ref{example33}. Specifically, as $k\to\sqrt{2L}-$,  the upper bound estimation for convergence time diverges, aligning with the numerical results in Example \ref{examplenew1}. Nevertheless, in Examples \ref{example22} and \ref{example33}, $\langle\bm{x},\bm{f}(\bm{x})\rangle$ is considerably smaller than $L\|\bm{x}\|^2$ in most of the time. Thus, our estimation is rough and leaves room for further improvement.
\end{remark}

\subsection{Synchronization}
Given that the synchronization problem can be viewed as a specific type of stabilization problem, it is possible to utilize the stochastic controller \eqref{2} for the purpose of achieving synchronization. More precisely, for the ordinary differential equations
\begin{equation}\label{4}
	{\rm d}{\bm x}_{t}=\bm f(\bm x_{t}){\rm d}t,
\end{equation}
we can rewrite~\eqref{1} as
\begin{equation}\label{3}
	{\rm d}\bm{y}_t=\bm{f}(\bm{y}_t){\rm d}t+\bm{u}(\bm{y}_t-\bm{x}_t){\rm d}B_t
\end{equation}
to achieve $\lim_{t\rightarrow+\infty}\bm z_{t}=0$ with $\bm z_{t}=\bm y_{t}-\bm x_{t}$. And we apply
	$$
\mathcal{E}_q\triangleq\int_0^\tau \|\bm{u}(\bm{y}_s-\bm{x}_s)\|^q{\rm d}s.
$$
as the energy consumption here.

\begin{example}\label{example_synchronization}
		 We consider the following Hindmarsh-Rose model \cite{b22, b32}:
		\begin{eqnarray*}
			{\rm d}{\bm x}_{t}&=&{\bm f}({\bm x}_{t}){\rm d}t, \\
			{\rm d}{\bm y}_{t}&=&{\bm f}({\bm y}_{t}){\rm d}t+{\bm u}({\bm y}_{t}-{\bm x}_{t}){\rm d}{B}_{t},
		\end{eqnarray*}
		where 
		\begin{equation*}
			\bm{f}(\bm x)=[
			x_{2}-x_{1}^{3}+3x_{1}^{2}-x_{3}+3, 1-5x_{1}^{2}-x_{2}, \epsilon(4x_{1}+6.4-x_{3})
			]^{\top}.
		\end{equation*}
	By defining $\bm z=\bm y-\bm x$, it can be observed that $\bm F(\bm z(t))=\bm f(\bm y)-\bm f(\bm x)=\bm f(\bm z+\bm x(t))-\bm f(\bm x(t))$ is not one-sided Lipschitzian with respect to $\bm z$. Nevertheless, the achievement of synchronization can be demonstrated through Figures~\ref{fig4} and \ref{fig5}. In Figure~\ref{fig5}, it can be observed that in the absence of coupling term, two systems exhibiting different initial values exhibit asynchronous spiking behavior. However, with the introduction of the coupling diffusion term $\bm u(\bm y_{t}-\bm x_{t})$, rapid synchronization between the two systems is achieved. Furthermore, the convergence time and energy consumption display a fluctuating decreasing trend as 
$k$ increases. It is likely that this behavior is influenced by errors stemming from the Euler-Maruyama scheme.

		\begin{figure}
		\includegraphics[width=0.45\textwidth]{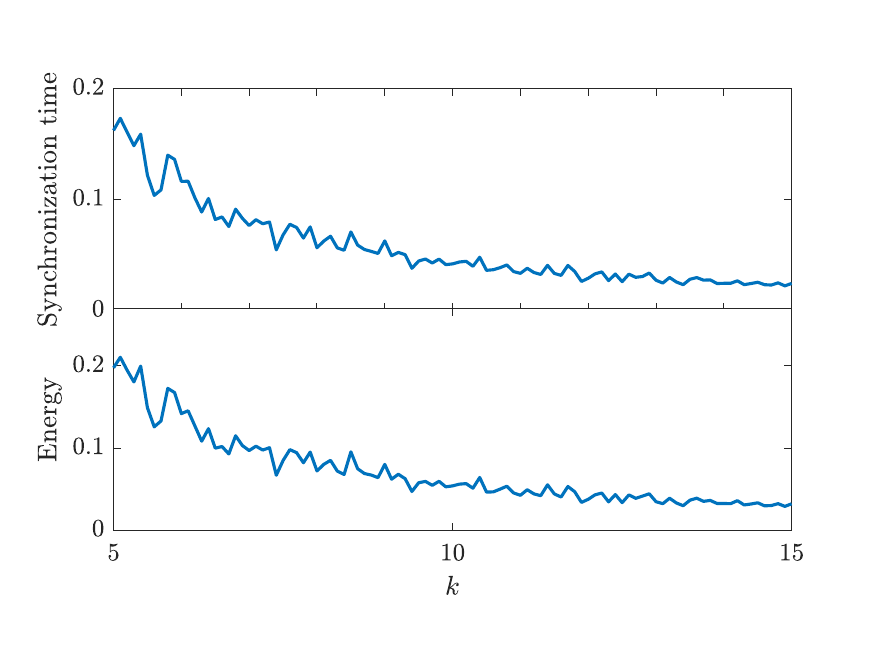}
		\caption{The convergence time and the required energy consumption for \eqref{example_synchronization}, respectively, vary with the increase of $k$. The parameters are $\epsilon=0.005$, $\alpha=0.5$, $q=0.1$, and the initial states are $\bm{x}_0=[0,0,1]^{\rm T}$, $\bm{x}_0=[0,0,2]^{\rm T}$. The solid (blue) curves indicate the numerical evaluation of $\tau$ and $\mathcal{E}_q$, which is counted as the average of 100 random realizations of the system \eqref{example_synchronization}. }\label{fig4}
	\end{figure}
		\begin{figure}
	\includegraphics[width=0.45\textwidth]{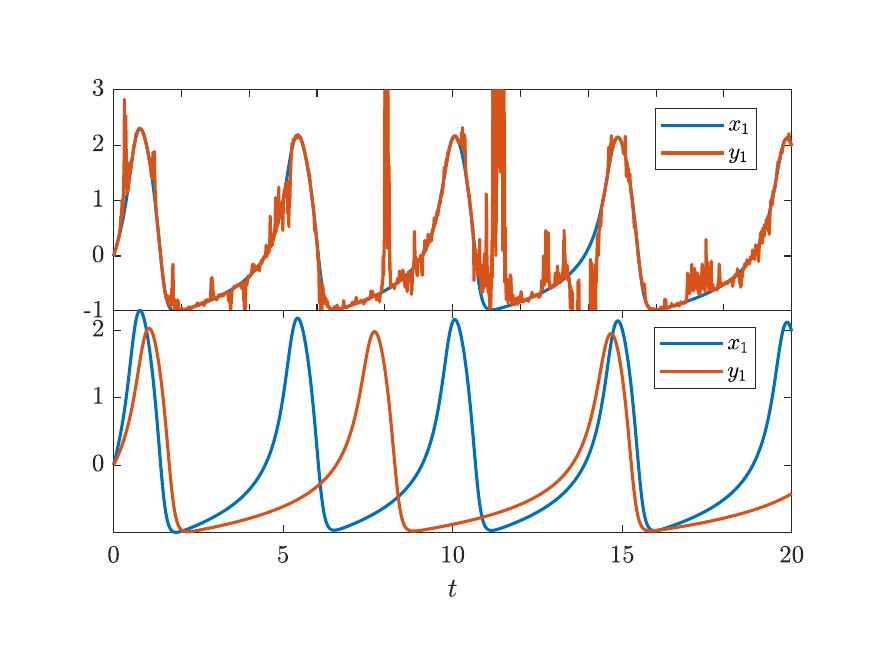}
	\caption{The dynamics of membrane potential $x_{1}$ and $y_{1}$ with (upper) and without (lower) the coupling diffusion term, respectively. The parameters are $\epsilon=0.005$, $\alpha=0.5$, $q=0.1$, and the initial states are $\bm{x}_0=[0,0,1]^{\rm T}$, $\bm{x}_0=[0,0,2]^{\rm T}$. Here, we only show one numerical realization of the the system \eqref{example_synchronization}.}\label{fig5}
\end{figure}
	
\end{example}

\section{Comparison with Deterministic Scheme}\label{comparison}
In this section, we compare our stochastic control scheme with the corresponding deterministic scheme. That is the substitution of~\eqref{1} as follows.
\begin{equation}\label{deterministic}
	{\rm d}\bm{x}_t=\bm{f}(\bm{x}_t){\rm d}t-\bm{u}(\bm{x}_t){\rm d}t,
\end{equation}
where $\bm u$ is also set in the form of~\eqref{2}. This scheme has been rigorously proved to be feasible for $k>L$, where $L$ represents the one-sided Lipschitz constant. When $L>2$, we observe $\sqrt{2L}<L$, indicating that the demand for $k$ is more stringent analytically for the deterministic scheme. Furthermore, convergence time and energy consumption, which are discussed and compared in the following example, serve as additional criteria for determining the superiority of either scheme.
\begin{example}\label{example3}
	Consider the Lorenz system \cite{b23, b31}, which reads
	\begin{eqnarray*}
		{\rm d}{\bm x}_{t}&=&{\bm f}({\bm x}_{t}){\rm d}t,
	\end{eqnarray*}
	where 
	\begin{equation*}
		\bm{f}(\bm x)=[
		\sigma x_{2}-\sigma x_{1},
		\rho x_{1}-x_{3}x_{1}-x_{2},
		x_{1}x_{2}-\beta x_{3} 
		]^{\top}.
	\end{equation*}
	
	Now, we have 
	\begin{eqnarray*}
		\langle{\bm x}, \bm f(\bm x)\rangle&=& (\sigma+\rho)x_{1}x_{2}-\sigma x_{1}^{2}- x_{2}^{2}-\beta x_{3}^{2}\\
		&\leq & \frac{\sigma+\rho}{2}(x_{1}^2 +x_{2}^2)\\
		&\leq &\frac{\sigma+\rho}{2} \vert \bm x \vert^{2}. 
	\end{eqnarray*}
	So the Lipschitzian constant $L$ can be taken as  $(\sigma+\beta)/2$ here. And if we take $k>\sqrt{\sigma+\rho}$ and $k>(\sigma+\beta)/2$, the Lorenz system can be stabilized in a finite time duration for the stochastic and deterministic scheme, respectively.  Figure~\ref{fig3} indicates the numerical simulation of convergence time and energy consumption for two different schemes. 
 
 It is clear that the convergence time of stochastic scheme is strictly smaller than the deterministic one. Besides, except for some specific values of $k$, the stochastic scheme consumes less energy than the deterministic one. This implies that the stochastic closed-loop controllers are generally better than the deterministic ones. But we also note that our scheme shows a large fluctuation for each trial, and it is greatly influenced by the time interval $\Delta t$ in Euler-Maruyama scheme (see \cite{b29, b30}). So there is much room for the improvement of robustness in our control scheme.
	
	\begin{figure}
		\includegraphics[width=0.45\textwidth]{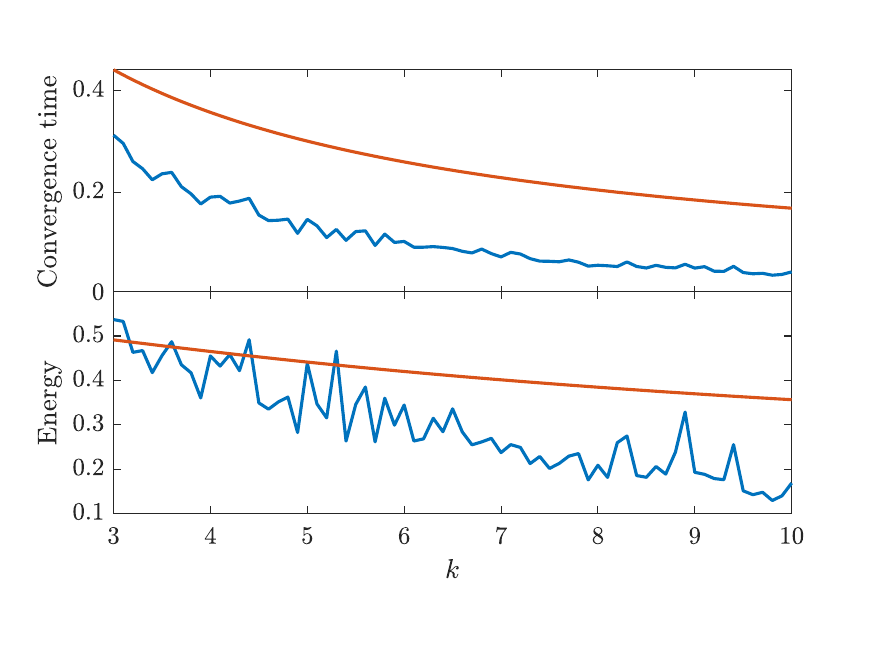}
		\caption{The convergence time and the required energy consumption for \eqref{example3}, vary respectively, with the increase of $k$. The parameters are $\sigma=6$, $\rho=10$, $\beta=3$, $\alpha=0.5$, $q=0.5$, and the initial states are $\bm{x}_0=[0,0,1]^{\rm T}$. The blue and red curves indicate the numerical evaluation of $\tau$ and $\mathcal{E}_q$ for the stochastic and deterministic control schemes, respectively, which is counted as the average of 200 random realizations of system \eqref{example3}. }\label{fig3}
	\end{figure}
	
\end{example}

\begin{example}\label{example5}
We now consider the following two-dimension systems:
\begin{equation*}
\left\{\begin{matrix}
\frac{\mathrm{d} x_{1}}{\mathrm{d} t}=\cos x_{2},   \\
\frac{\mathrm{d} x_{2}}{\mathrm{d} t}=\sin x_{1}, 
\end{matrix}\right.
\end{equation*}
which starts from the state $[0,1]^{\top}$. It obviously satisfies the global Lipschitzian condition. For another system with the same drift term but starts from different initial value $\bm y_{t}$, we apply the scheme ${\rm d}{\bm y}_{t}={\bm f}({\bm y}_{t}){\rm d}t+{\bm u}({\bm y}_{t}-{\bm x}_{t}){\rm d}{B}_{t}$ and ${\rm d}{\bm y}_{t}={\bm f}({\bm y}_{t}){\rm d}t-{\bm u}({\bm y}_{t}-{\bm x}_{t}){\rm d}{t}$, respectively, to achieve synchronization. We compare their synchronization time and the energy consumption in Figure~\ref{synchronization_comparison}.

\begin{figure}
		\includegraphics[width=0.45\textwidth]{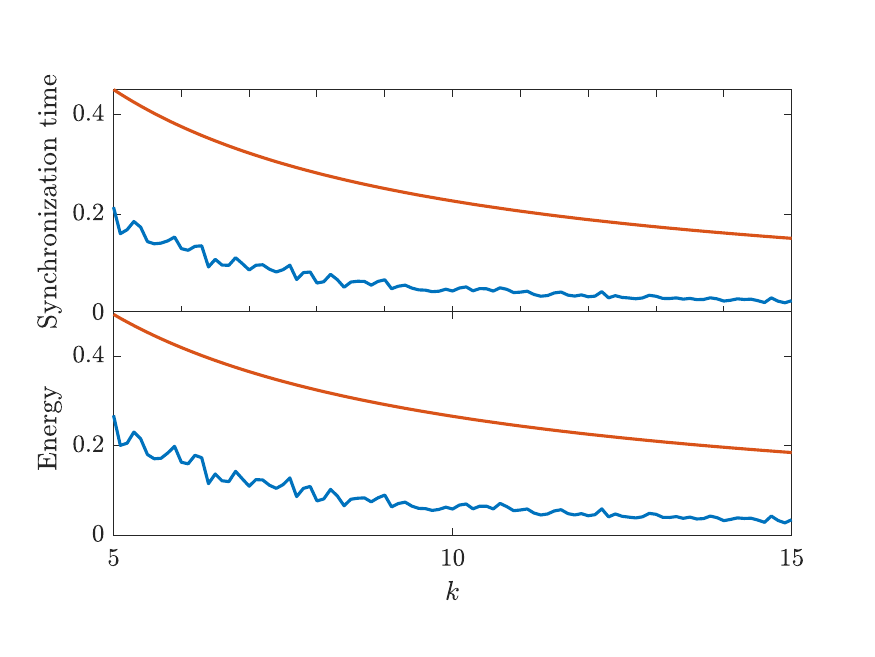}
		\caption{The synchronization time and the required energy consumption for \eqref{example5}, vary respectively, with the increase of $k$. The parameters are  $\alpha=0.5$, $q=0.1$, and the initial states are $\bm{x}_0=[0,1]^{\rm T}$ and $\bm{y}_0=[1,0]^{\rm T}$. The blue and red curves indicate the numerical evaluation of $\tau$ and $\mathcal{E}_q$ for the stochastic and deterministic control schemes, respectively, which is counted as the average of 100 random realizations of system \eqref{example5}. }\label{synchronization_comparison}
	\end{figure}

 The stochastic scheme here still shows great progress compared to the deterministic scheme. It always shows a shorter transition to synchronization as well as less energy consumption. This probably comes from some large value of $\Delta B_{t}$, which can synchronize two systems drastically in a short period. And a shorter synchronization time often means a less energy dissipation. Thus, the volatility here contributes greatly to the close-loop control.
\end{example}

\section{Analytical Validation}\label{proof}

The following generalized It$\hat{o}$'s formula, which can be found in \cite[Problem 7.3]{b1}, will be used in this section.
\begin{lemma}\label{lemma1}
	Let 
	$$
	X_t=X_0+M_t+V_t, 0\leq t<+\infty
	$$
	be a continuous semimartingale, where $M=\{M_t, \mathscr{F}_t; 0\leq t<+\infty\}$ is a local martingale, and $V=\{V_t,\mathscr{F}_t;0\leq t<+\infty\}$ is the difference of \textbf{continuous}, nondecreasing adapted process with $V_0=0$, and $\{\mathscr{F}_t\}$ satisfies the usual conditions. If $f:\mathbb{R}\to\mathbb{R}$ be a function whose derivative is \textbf{absolutely continuous}. Then $f''$ exists Lebesgue-almost everywhere, and we have the generalized It$\hat{o}$'s formula:
	$$
	\begin{aligned}
		f(X_t)=&f(X_0)+\int_0^t f'(X_s){\rm d}M_s+\int_0^t f'(X_s){\rm d}V_s\\
		&+\dfrac{1}{2}\int_0^t f''(X_s){\rm d}\langle M \rangle_s~~{\rm a.s.},
	\end{aligned}
	$$
	where $\langle M \rangle$ denotes by the quadratic process of $M$. 
\end{lemma}
\subsection{Estimation of Control Time}\label{time}

{\em Proof of Theorem \ref{theorem1}:} 
 We introduce the following function,
$$
V_p(x)=\left\{
\begin{aligned}
	&p\log x+1, x>1,
	\\&x^p, 0\leq x<1.
\end{aligned}\right.
$$
We note that $V$ is a continuous function defined on $[0,+\infty)$ with absolutely continuous first derivative on $(0,+\infty)$. Here, $p\in(0,1)$ is a parameter to be identified later. Therefore, an application of the generalized It$\hat{o}$'s formula (Lemma \ref{lemma1}) yields that
$$
V_p(\|\bm{x}_t\|)-V_p(\|\bm{x}_0\|)=\int_0^t \mathscr{L}V_p(\bm{x}_s){\rm d}s+\int_0^t \mathscr{D}V_p(\bm{x}_s){\rm d}B_s,
$$
where 
$$
\mathscr{L}V_p(\bm{x})\triangleq\left\{
\begin{aligned}
	&p\Bigg[\dfrac{\mathscr{G}(\bm{x})}{\|\bm{x}\|^2}-\dfrac{k^2}{2}\Bigg],\|\bm{x}\|\geq 1,\\
	&p\|\bm{x}\|^{p-2}\mathscr{G}(\bm{x})+\dfrac{p(p-1)k^2}{2}\|\bm{x}\|^{p-2+2\alpha},\\
	&\|\bm{x}\|<1.
\end{aligned}\right.
$$

$$
\hspace{-3cm}\mathscr{D}V_p(\bm{x})\triangleq\left\{
\begin{aligned}
	&kp,\|\bm{x}\|\geq 1,\\
	&kp\|\bm{x}\|^{\alpha+p-1},\|\bm{x}\|<1.
\end{aligned}\right.
$$
Here, $\mathscr{G}(\bm{x})\triangleq\langle\bm{x},\bm{f}(\bm{x}) \rangle\leq L\|\bm{x}\|^2$. Choosing $0<p\leq\min\left\{1,2-2\alpha\right\}$ and using $k>\sqrt{2L}$ yield following estimation:
\begin{equation}\label{proofjia2}
\mathscr{L}V_p(\bm{x})\leq\left\{
\begin{aligned}
	&-p\left(\dfrac{k^2}{2}-L\right)\triangleq H_1(p),\|\bm{x}\|\geq 1,\\
	& H_2(p),\|\bm{x}\|<1,
\end{aligned}\right.
\end{equation}
where $H_2(p)\triangleq pL+\dfrac{p(p-1)k^2}{2}$. Denote by a stopping time $\sigma\triangleq \inf \Bigg\{t>0\Bigg|\|\bm{x}_t\|=1\Bigg\}$.  Apparently, we obtain that $\sigma\leq\tau$. Then we give the estimation for $\tau$ on two disjoint measurable sets $\bm{1}_{\|\bm{x}_0\|>1}$ and $\bm{1}_{\|\bm{x}_0\|\leq 1}$ seperately.

I: The estimation for the set $\bm{1}_{\|\bm{x}_0\|>1}$.
An application of the generalized It$\hat{o}$'s formula ( Lemma \ref{lemma1}) for $V_p(\|\bm{x}_t\|)$ on the time interval $[0,t\wedge\sigma]$ yields that
$$
\begin{aligned}
	V_p(\|\bm{x}_{t\wedge\sigma}\|)-V_p(\|\bm{x}_0\|)=\int_0^{t\wedge\sigma} \mathscr{L}V_p(\bm{x}_s){\rm d}s\\
 +\int_0^{t\wedge\sigma} \mathscr{D}V_p(\bm{x}_s){\rm d}B_s.
\end{aligned}
$$
Taking the expectation of both sides on the measurable set $\bm{1}_{\|\bm{x}_0\|>1}$ yields that
$$
\begin{aligned}
\mathbb{P}(\|\bm{x}_0\|>1)\leq&\mathbb{E}V_p(\|\bm{x}_{t\wedge\sigma}\|;\bm{1}_{\|\bm{x}_0\|>1})\\
\leq&\mathbb{E}V_p(\|\bm{x}_0\|;\bm{1}_{\|\bm{x}_0\|>1})\\
&+H_1(p)\mathbb{E}(t\wedge\sigma;\bm{1}_{\|\bm{x}_0\|>1}),
\end{aligned}
$$
which further indicates that
$$
\mathbb{E}(t\wedge\sigma;\bm{1}_{\|\bm{x}_0\|>1})\leq \dfrac{2\mathbb{E}(\log\|\bm{x}_0\|;\bm{1}_{\|\bm{x}_0\|>1})}{{k^2}-2L}.
$$
Letting $t\to+\infty$ leads to that 
$$
\mathbb{E}(\sigma;\bm{1}_{\|\bm{x}_0\|>1})\leq \dfrac{2\mathbb{E}(\log\|\bm{x}_0\|;\bm{1}_{\|\bm{x}_0\|>1})}{{k^2}-2L}
$$
and $\sigma<+\infty$ almost surely.

Moreover, an application of the  generalized It$\hat{o}$'s formula (Lemma \ref{lemma1}) for $V_p(\|\bm{x}_t\|)$ on the time interval $[\sigma,(\sigma+t)\wedge\tau]$ yields that
\begin{equation}\label{proof1}
\begin{aligned}
	V_p(\|\bm{x}_{(\sigma+t)\wedge\tau}\|)-V_p(\|\bm{x}_{\sigma}\|)=&\int_\sigma^{(\sigma+t)\wedge\tau} \mathscr{L}V_p(\bm{x}_s){\rm d}s\\&+\int_\sigma^{(\sigma+t)\wedge\tau} \mathscr{D}V_p(\bm{x}_s){\rm d}B_s.
\end{aligned}
\end{equation}
Note that even for $t>\sigma$, the trajectory may be outside the unit ball. However, we have the estimation $H_1(p)<H_2(p)$, which implies that $ 
\mathscr{L}V_p(\bm{x}_s)\leq H_2(p)~~{\rm a.s.}.$  Taking the expectation of both sides yields that
$$
\begin{aligned}
0&\leq \mathbb{E}V_p(\|\bm{x}_{(\sigma+t)\wedge\tau}\| 
;\bm{1}_{\|\bm{x}_0\|>1})\\
&\leq \mathbb{P}(\|\bm{x}_0\|>1)+H_2(p)\mathbb{E}\left[(\sigma+t)\wedge\tau-\sigma;\bm{1}_{\|\bm{x}_0\|>1}\right].
\end{aligned}
$$
When $H_2(p)<0$, we obtain that
$$
\mathbb{E}\left[(\sigma+t)\wedge\tau-\sigma 
;\bm{1}_{\|\bm{x}_0\|>1}\right]\leq -\dfrac{\mathbb{P}(\|\bm{x}_0\|>1)}{H_2(p)}.
$$
Letting $t\to+\infty$ leads to that 
$$
\mathbb{E}\left[\tau-\sigma;\bm{1}_{\|\bm{x}_0\|>1}\right]\leq -\dfrac{\mathbb{P}(\|\bm{x}_0\|>1)}{H_2(p)}
$$
and $\tau<+\infty$ almost surely. We obtain that 
$$ 
\mathbb{E}(\tau;\bm{1}_{\|\bm{x}_0\|>1})\leq \dfrac{2\mathbb{E}(\log\|\bm{x}_0\|;\bm{1}_{\|\bm{x}_0\|>1})}{{k^2}-2L}-\dfrac{\mathbb{P}(\|\bm{x}_0\|>1)}{H_2(p)}.$$
\textcolor{black}{II: The estimation for the set $\bm{1}_{\|\bm{x}_0\|\leq 1}.$ On this set we have $\sigma=0$. Eq. \eqref{proof1}, constrained within the set $\bm{1}_{\|\bm{x}_0\|\leq 1}$, yields that
\begin{equation}\label{proofjia1}
\begin{aligned}
	&V_p(\|\bm{x}_{t\wedge\tau}\|; \bm{1}_{\|\bm{x}_0\|\leq 1})-V_p(\|\bm{x}_{0}\|; \bm{1}_{\|\bm{x}_0\|\leq 1})\\&=\int_0^{t\wedge\tau} \mathscr{L}V_p(\bm{x}_s)\bm{1}_{\|\bm{x}_0\|\leq 1}{\rm d}s\\&+\int_0^{t\wedge\tau} \mathscr{D}V_p(\bm{x}_s)\bm{1}_{\|\bm{x}_0\|\leq 1}{\rm d}B_s.
\end{aligned}
\end{equation}
Note that the estimation \eqref{proofjia2} and $H_1(p)\leq H_2(p)$ still apply, we have  $\mathscr{L}V_p(\bm{x}_s)\leq H_2(p)~~{\rm a.s.}$ even $\bm{x}_t$ may be outside the unit ball. By taking expectation on Eq. \eqref{proofjia1}, we have
$$
\begin{aligned}
&\hspace{0.5cm}-\mathbb{E}V_p(\|\bm{x}_{0}\|; \bm{1}_{\|\bm{x}_0\|\leq 1}) \\
&\leq \mathbb{E}[V_p(\|\bm{x}_{t\wedge\tau}\|; \bm{1}_{\|\bm{x}_0\|\leq 1})-V_p(\|\bm{x}_{0}\|; \bm{1}_{\|\bm{x}_0\|\leq 1})] \\
&\leq H_2(p)\mathbb{E}(t\wedge\tau; \bm{1}_{\|\bm{x}_0\|\leq 1}).
\end{aligned}
$$
When $H_2(p)<0$, we obtain that
$$
\begin{aligned}
\mathbb{E}(t\wedge\tau; \bm{1}_{\|\bm{x}_0\|\leq 1})&\leq -\dfrac{\mathbb{E}V_p(\|\bm{x}_{0}\|; \bm{1}_{\|\bm{x}_0\|\leq 1})}{H_2(p)}\\
&=-\dfrac{\mathbb{E}(\|\bm{x}_0\|^p;\bm{1}_{\|\bm{x}_0\|\leq 1})}{H_2(p)}.
\end{aligned}
$$
Letting $t\to+\infty$ yields that
$$
\mathbb{E}(\tau;\bm{1}_{\|\bm{x}_0\|\leq 1})\leq -\dfrac{\mathbb{E}(\|\bm{x}_0\|^p;\bm{1}_{\|\bm{x}_0\|\leq 1})}{H_2(p)}.
$$}
Consequently, we have 
$$
\begin{aligned}
\mathbb{E}\tau\leq & \dfrac{2\mathbb{E}(\log\|\bm{x}_0\|;\bm{1}_{\|\bm{x}_0\|>1})}{{k^2}-2L}-\dfrac{\mathbb{P}(\|\bm{x}_0\|>1)}{H_2(p)}\\
&-\dfrac{\mathbb{E}(\|\bm{x}_0\|^p;\bm{1}_{\|\bm{x}_0\|\leq 1})}{H_2(p)}.
\end{aligned}
$$

Now we come to investigate $H_2(p)$. Since $H_2'(p)=pk^2-\left(\dfrac{k^2}{2}-L\right)$, $H_2(p)$ decreases on the interval $(0,p^*)$ and increases on the interval $(p^*,+\infty)$, where $p^*\triangleq\dfrac{1}{2}-\dfrac{L}{k^2}$. Therefore, we obtain that 
$$
\min_{0<p\leq \min\{1,2-2\alpha\}}H_2(p)=\left\{\begin{aligned}
	-\dfrac{{(kp^*)}^2}{2}, p^*<2-2\alpha, \\
	H_2(2-2\alpha), p^*\geq 2-2\alpha.\end{aligned}\right.
$$
If we choose $p=p^*$ when $\alpha<\dfrac{3}{4}+\dfrac{L}{2k^2}$ and $p=2-2\alpha$ when $\alpha>\dfrac{3}{4}+\dfrac{L}{2k^2}$, we get the following estimations.
$$
\mathbb{E}\tau\leq T_{f}^{\rm Sup}\triangleq\left\{
\begin{aligned}
	&\dfrac{2\mathbb{E}(\log\|\bm{x}_0\|;\bm{1}_{\|\bm{x}_0\|>1})}{{k^2}-2L}\\
 &+\dfrac{8k^2}{(k^2-2L)^2}[\mathbb{P}(\|\bm{x}_0\|>1)
\\
&+\mathbb{E}(\|\bm{x}_0\|^{p^*};\bm{1}_{\|\bm{x}_0\|\leq 1})], \
 \alpha<\dfrac{3}{4}+\dfrac{L}{2k^2},\\
	&\dfrac{2\mathbb{E}(\log\|\bm{x}_0\|;\bm{1}_{\|\bm{x}_0\|>1})}{{k^2}-2L}\\
 &+\dfrac{[\mathbb{P}(\|\bm{x}_0\|>1)
+\mathbb{E}(\|\bm{x}_0\|^{2-2\alpha};\bm{1}_{\|\bm{x}_0\|\leq 1})]}{(1-\alpha)[(2\alpha-1)k^2-2L]},\\
	&\alpha\geq \dfrac{3}{4}+\dfrac{L}{2k^2}.
\end{aligned}\right.
$$

\hfill{$\blacksquare$}

\subsection{Estimation of Control Energy}\label{energy}

{\em Proof of Theorem \ref{theorem2}:} 
The energy consumption of the control is given by

\begin{eqnarray*}
\mathcal{E}_q&=&\int_0^\tau \|\bm{u}_s\|^q{\rm d}s\leq k^q\int_0^\sigma \|\bm{x}_s\|^q{\rm d}s\\
&&+k^q\int_\sigma^\tau (1 \vee \|\bm{x}_s\|^{q }){\rm d}s\\
&\leq &k^q\int_0^\sigma \|\bm{x}_s\|^q{\rm d}s+k^q\int_\sigma^\tau (1 + \|\bm{x}_s\|^{q }){\rm d}s
\end{eqnarray*}
Denote by $W_q(x)\triangleq x^q$. An application of the generalized It$\hat{o}$'s formula (Lemma \ref{lemma1}) yields that
\begin{eqnarray*}
W_q(\|\bm{x}_t\|)-W_q(\|\bm{x}_0\|)&=&\int_0^t \mathscr{L}W_q(\bm{x}_s){\rm d}s\\
&&+\int_0^t\mathscr{D}W_q(\bm{x}_s){\rm d}B_s,
\end{eqnarray*}
where

\begin{equation}
	\mathscr{L}W_q(\bm{x})=\left\{
	\begin{aligned}
		&q\|\bm{x}\|^{q-2}\mathscr{G}(\bm{x})+\dfrac{q(q-1)k^2}{2}\|\bm{x}\|^q, \ \|\bm{x}\|\geq 1,\\
	&\|\bm{x}\|^{q-2}\mathscr{G}(\bm{x})+\dfrac{q(q-1)k^2}{2}\|\bm{x}\|^{q-2+2\alpha}, \\
  &\|\bm{x}\|<1.
	\end{aligned}\right.
\end{equation}
and
$$
\hspace{-3cm}\mathscr{D}W_q(\bm{x})\triangleq\left\{
\begin{aligned}
	&kq\|\bm{x}\|^{q},\|\bm{x}\|\geq 1,\\
	&kq\|\bm{x}\|^{\alpha+q-1},\|\bm{x}\|<1.
\end{aligned}\right.
$$

Choosing $0<q<\min\Bigg\{2-2\alpha,1-\dfrac{2L}{k^2}\Bigg\}$ yields following estimations

\begin{equation}
	\mathscr{L}W_q(\bm{x})\leq\left\{
	\begin{aligned}
		& H_2(q)\|\bm{x}\|^q    ,\|\bm{x}\|\geq 1,\\
		&H_2( q),\|\bm{x}\|<1.
	\end{aligned}\right.
\end{equation}
I: We first do the estimation on the set $\bm{1}_{\|\bm{x}_0\|>1}$. Thus, on the time interval $[0,t\wedge\sigma]$, we have 
\begin{eqnarray*}
	W_q(\|\bm{x}_{t\wedge\sigma}\|)&=&W_q(\|\bm{x}_0\|)+\int_0^{t\wedge\sigma} \mathscr{L}W_q(\bm{x}_s){\rm d}s\\
	&&+\int_0^{t\wedge\sigma}\mathscr{D}W_q(\bm{x}_s){\rm d}B_s.
\end{eqnarray*}

Taking the expectation of both sides on $\bm 1_{\|\bm x_{0}\|>1}$ leads to that
$$
\begin{aligned}
	\mathbb{E}\bigg[\|\bm{x}_t\|^q\bm{1}_{t<\sigma,\|\bm x_{0}\|>1}\bigg]&\leq\mathbb{E}\left [W_q(\|\bm{x}_{t\wedge\sigma}\|);\bm {1}_{\|\bm x_{0}\|>1}\right ]\\&\leq \mathbb{E} \left [ W_q(\|\bm{x}_0\|);\bm {1}_{\|\bm x_{0}\|>1} \right ]\\
 &+H_2(q)\int_0^t \mathbb{E}\bigg[\|\bm{x}_s\|^q\bm{1}_{s<\sigma,\|\bm x_{0}\|>1}\bigg]{\rm d}s,
\end{aligned}
$$
this, together with Gronwell inequality, leads to that
$$
\begin{aligned}
&\mathbb{E}\bigg[\|\bm{x}_t\|^q\bm{1}_{t<\sigma,\|\bm x_{0}\|>1}\bigg]\\
\leq& \mathbb{E} \left [ W_q(\|\bm{x}_0\|);\bm {1}_{\|\bm x_{0}\|>1}\right ]\exp\bigg[H_2(q)t\bigg].
\end{aligned}
$$
Note that $H_2( q)<0$, we obtain 
$$
\begin{aligned}
	&\mathbb{E}\Bigg [\int_0^\sigma \|\bm{x}_s\|^q{\rm d}s ;\|\bm x_{0}\|>1 \Bigg ]\\
	=&\int_0^{+\infty}\mathbb{E}\bigg[\|\bm{x}_t\|^q\bm{1}_{t<\sigma,\|\bm x_{0}\|>1}\bigg]{\rm d}t \\
	\leq& \mathbb{E} \left [ W_q(\|\bm{x}_0\|);\bm {1}_{\|\bm x_{0}\|>1} \right ]\int_0^{+\infty}\exp\bigg[H_2(q)t\bigg]{\rm d}t\\
 \leq& -\dfrac{\mathbb{E} \left [ W_q(\|\bm{x}_0\|);\bm {1}_{\|\bm x_{0}\|>1} \right ]}{H_2(q)}.
\end{aligned}
$$
Furthermore, we have
$$
\begin{aligned}
	W_q(\|\bm{x}_{(t+\sigma)\wedge\tau}\|)=&W_q(\|\bm{x}_\sigma\|)+\int_\sigma^{(t+\sigma)\wedge\tau} \mathscr{L}W_q(\bm{x}_s){\rm d}s\\&+\int_\sigma^{{(t+\sigma)\wedge\tau}}\mathscr{D}W_q(\bm{x}_s){\rm d}B_s,
\end{aligned}
$$
which implies that
\begin{eqnarray*}
&&\mathbb{E}\bigg[\|\bm{x}_{t+\sigma}\|^q\bm{1}_{t<\tau-\sigma,\|\bm x_{0}\|>1}\bigg]\\&=& \mathbb{E}\bigg[W_q(\|\bm{x}_{(t+\sigma)\wedge\tau}\|)\bm{1}_{t<\tau-\sigma,\|\bm x_{0}\|>1}\bigg]\\
&\leq& \mathbb{E}[\bm{1}_{t<\tau-\sigma,\|\bm x_{0}\|>1}].
\end{eqnarray*}
Therefore, we obtain that
$$
\begin{aligned}
	&\mathbb{E}\left [\int_\sigma^\tau \|\bm{x}_s\|^q{\rm d}s;\|\bm x_{0}\|>1 \right ]\\&=\int_0^{\tau-\sigma}\mathbb{E}\bigg[\|\bm{x}_{t+\sigma}\|^q\bm{1}_{0<t<\tau-\sigma}\bigg]{\rm d}t\\&\leq \int_0^{+\infty}\mathbb{E}\bigg[\bm{1}_{0<t<\tau-\sigma,\|\bm x_{0}\|>1}\bigg]{\rm d}t\\&=\mathbb{E}\bigg[\int_0^{+\infty}\bm{1}_{0<t<\tau-\sigma}
	{\rm d}t; \|\bm x_{0}\|>1\bigg]\\
	&=\mathbb{E}\bigg[\tau-\sigma ;\|\bm x_{0}\|>1\bigg]\leq -\dfrac{\mathbb{P}(\|\bm{x}_0\|>1)}{H_2(q)}.
\end{aligned}
$$
This results in that

\begin{eqnarray*}
&&\mathbb{E}[\mathcal{E}_q;\|\bm{x}_0\|>1]\\&\leq& k^q \mathbb{E}\Bigg [\int_0^\sigma \|\bm{x}_s\|^q{\rm d}s;\|\bm{x}_0\|>1\Bigg ]\\
&&+k^q \mathbb{E}\Bigg [\int_\sigma^\tau (1+\|\bm{x}_s\|^q){\rm d}s;\|\bm{x}_0\|>1\Bigg ]\\
&\leq & -k^q\dfrac{\mathbb{E}[\|\bm{x}_0\|^q+2;\|\bm{x}_0\|>1]}{H_2(q)}.
\end{eqnarray*}

II: Secondly, we do estimations on the set  $\bm{1}_{\|\bm{x}_0\|\leq 1}$.
A similar argument leads to that

$$
\mathbb{E}[\mathcal{E}_q;\|\bm{x}_0\|\leq 1]\leq -\dfrac{2k^q \mathbb{E}[\|\bm{x}_0\|^q; \|\bm{x}_0\|\leq 1)]}{H_2(q)}.
$$

Above all, we have the estimation
$$
\begin{aligned}
\mathbb{E}\mathcal{E}_q\leq\mathcal{E}_q^{\rm Sup}&\triangleq-\dfrac{k^q}{H_2 (q)}[\mathbb{E}(\|\bm x_0\|^q )+2\mathbb{P}(\|\bm x_0\|>1)\\
&+\mathbb{E}(\|\bm x_0\|^q ;\|\bm x_0\|\leq 1)]
 \end{aligned}
	$$
 This completes the proof of the whole theorem.  \hfill{$\blacksquare$}

\textcolor{black}{\begin{remark}\label{remarkjia1}
The conclusion that Eq. \eqref{1} is finite-time stable ($\mathbb{P}(\tau<+\infty)=1$) is obtained by using the result of Theorem \ref{theorem1}, which states that $\mathbb{E}\tau<+\infty$. In most literatures (for example, see \cite{b50,b52}), the definition for finite-time stable is $\mathbb{P}(\tau<+\infty)=1$ (refer to Definition \ref{def1}) rather than $\mathbb{E}\tau<+\infty$. This definition is commonly used because it is less restrictive and easier to verify in certain cases.  For example, we consider a time-varying one-dimensional SDE (refer to \cite[Example 3.1]{b52})
\begin{equation}\label{lizijia1}
{\rm d}x=(-\dfrac{a}{2(1+t)}x^{\frac{1}{3}}+\dfrac{b}{1+t^2}x){\rm d}t+\dfrac{d}{\sqrt{1+t}}x^{\frac{2}{3}}{\rm d}B_t
\end{equation}
where the parameters $a\geq d^2, b\geq 0$. Authors in \cite{b52} show that Eq. \eqref{lizijia1} is finite-time stable by using \cite[Corollary 3.2]{b52} and obtaining that
$
\mathbb{E}\beta(\tau)<+\infty,
$
where $\beta(t)=\int_0^t c(s)({\rm e}^{-\int_0^s l(u){\rm d}u})^{1-\gamma}{\rm d}s, c(t)=\dfrac{a}{2(1+t)}, \gamma=\frac{1}{3}, l(t)=\dfrac{b}{1+t^2}$. However, it is challenging to demonstrate that $\mathbb{E}\tau<+\infty$ in this particular case.
\end{remark}}

\begin{remark}\label{remark2}
We need to emphasize that even for $t>\sigma$, the trajectory may be outside the unit circle because it is a stochastic system. Therefore, we cannot seperate the time interval into two parts and do the estimation using different Lyapunov functions seperately.  This is totally different from the past work on the deterministic systems (See \cite{b15}). This is reason why we construct a piecewise smooth Lyapnounv function and use the generalized It$\hat{o}$'s formula (Lemma \ref{lemma1}), which is the major difficulty for the proof.
\end{remark}

\section{Concluding Remarks}\label{conclusion}
In this paper, we propose a finite-time, closed-loop stochastic feedback controller. The usefulness of the scheme in chaos control and network synchronization is illustrated not only by analytical validations but also by numerically controlling several representative stochastic dynamical models. Additionally, we investigate convergence time and energy consumption as well as its dependence on parameters. And we also numerically validate the our estimation of the upper bound.

As for the future research directions, the estimation of convergence time and energy consumption is much larger than the exact value by the simulation. Moreover, our findings suggest that under certain conditions of $\alpha$,  the values of $T_{f}^{\rm Sup},\mathcal{E}_q^{\rm Sup}$ remains a constant when $\alpha$ satisfies particular conditions.Thus, there is a need for more accurate estimation techniques to address this issue.

In addition, there are practical implementation challenges associated with our proposed controller. One particular challenge is the reliance on instant-time information, which may not always be available in practical applications. To address this issue, event-triggered schemes have been developed for stabilizing deterministic dynamical systems (See \cite{b64, b65, b66, b67}). Additionally, schemes incorporating time delay have also been explored (See \cite{b70}). Furthermore, in practical applications, uncertainties such as unknown parameters, unmodeled dynamics, uncertain functions, and disturbances are inevitable. For stabilizing deterministic systems, several finite-time adaptive control designs have been established to address uncertainty (See \cite{b63, b68}). It would be worthwhile to investigate how to extend these results to stochastic controllers in future research.

\section*{Data availability statement}
The data that supports the findings of this study are computationally generated and available upon request.

\section*{Acknowledgments}
The authors thank Luan Yang for plotting Figure 1. S.Z.
is supported by the National Natural Science Foundation
of China (No. 11925103) and by the STCSM (Nos.
22JC1402500, 22JC1401402, and 2021SHZDZX0103). He is also supported by the Research on the Assessment and Prevention Techniques for Telecom Network Fraud Victims Based on Big Data (No. 21DZ1201402).

\bibliographystyle{cas-model2-names}

\bibliography{bibliography}

\end{document}